\documentclass[12pt,a4paper]{amsart}
 \usepackage[T1]{fontenc}
\usepackage[top=35mm, bottom=35mm, left=30mm, right=30mm]{geometry}    

\usepackage{amssymb}
\usepackage{mathtools}
\usepackage{verbatim,enumerate}
\usepackage{xcolor}

\usepackage[looser]{newpxtext}
\usepackage[smallerops]{newpxmath}

\usepackage{fancyhdr}
\makeatletter
\def\@@and{\MakeLowercase{and}}
\makeatother

\usepackage[hidelinks,pagebackref]{hyperref}

\usepackage[nobysame]{amsrefs} 

\theoremstyle{definition}
\newtheorem{defn}{Definition}[section]

\newtheorem{question}[defn]{Question}
\newtheorem{rem}[defn]{Remark}

\theoremstyle{plain}
\newtheorem{thm}[defn]{Theorem}
\newtheorem{lem}[defn]{Lemma}
\newtheorem{prop}[defn]{Proposition}
\newtheorem{coro}[defn]{Corollary}
 
\newcommand{\bbz}{\mathbb{Z}} 
\newcommand{\bbn}{\mathbb{N}}
\DeclareMathOperator{\diam}{diam}

\title[D\MakeLowercase{ense  uniform} L\MakeLowercase{i}-Y\MakeLowercase{orke chaos for linear operators}] 
{D\MakeLowercase{ense uniform} L\MakeLowercase{i}-Y\MakeLowercase{orke chaos for linear operators on a }B\MakeLowercase{anach space}}

\author[J. L\MakeLowercase{i}]{J\MakeLowercase{ian} Li}
\address[J. Li]{Department of Mathematics,
	Shantou University, Shantou, 515821, Guangdong, China}
\email{lijian09@mail.ustc.edu.cn}
\urladdr{https://orcid.org/0000-0002-8724-3050}

\author[X. W\MakeLowercase{ang}]{X\MakeLowercase{insheng} Wang}
\address[X. Wang]{Department of Mathematics,
	Shantou University, Shantou, 515821, Guangdong, China}
\email{wangxs@stu.edu.cn}
\urladdr{https://orcid.org/0000-0002-5287-8902}

	

\date{\today}

\begin{document}

\begin{abstract}
	This paper focuses on the dense uniform Li-Yorke chaos for linear operators on a Banach space. Some sufficient conditions and equivalent conditions are established under which the dynamical system is densely uniformly Li-Yorke chaotic. It is shown that there are plenty of densely uniformly Li-Yorke chaotic operators. For unilateral backward weighted shifts and  bilateral backward weighted shifts on $\ell^p$, it is shown that Li-Yorke chaos is equivalent to dense uniform Li-Yorke chaos.
\end{abstract}

\maketitle


\section{Introduction}
It is an interesting and meaningful question how complicated a system can be in the study of dynamical systems. Plenty of concepts and invariants have been introduced, giving quantitative or qualitative descriptions of the complexity for various systems from different points of view. Among these concepts, chaos is an important one. Note that there are multiple different types of definitions for chaos, however these definitions are not equivalent to each other. In \cite{LY1975} Li and Yorke introduced one kind of chaos for continuous maps on the interval, and showed that if a map has points with period three, then it is chaotic. In the following, we will give the definition of Li-Yorke chaos under a more general setting.

Let $(X, d)$ be a metric space. A continuous map $f \colon X \to X$ is called \emph{Li–Yorke chaotic} if there exists an uncountable subset $S\subset X$ such that for every pair $x,y\in S$ with $x\neq y$, we have
\[
	\liminf_{n\to\infty} d(f^n(x),f^n(y))=0\text{ and }
	\limsup_{n\to\infty} d(f^n(x),f^n(y))>0.
\]
Intuitively speaking, the distance between distinct points in $S$ is sometimes far and sometimes close driven by $f$.

Another popular definition of chaos is given by Devaney. A continuous map $f\colon X\to X$ is called \emph{Devaney chaotic} if it satisfies the following conditions:
\begin{enumerate}
	\item $f$ is transitive (see the paragraph below Corollary \ref{coro:Dk} for the definition);
	\item the set of periodic points of $f$ is dense in $X$;
	\item $f$ is sensitive dependence on initial conditions (see the paragraph above Lemma \ref{lem:T-sensitive-eq} for the definition).
\end{enumerate}
In \cite{HY2002} Huang and Ye proved that Devaney chaos implies Li-Yorke chaos, which is an application based on their deep investigation of the asymptotic cell and relation and  sufficient conditions under which the system is Li-Yorke chaotic. Generally, they proved that if $f$ is transitive and $X$ is infinite containing a periodic point, then there is an uncountable scrambled set  for $f$ (the definition of scrambled set will be given at the beginning of Section \ref{sec:DULYmetric}). Note that the \emph{diameter} of a nonempty subset $A$ of $X$ is defined by
\[
	\diam(A)=\sup_{x,y\in A}d(x,y).
\]
Furthermore, in \cite{M2004} Mai strengthened Huang and Ye's result, showing that if $f$ is transitive and has a periodic point with period $p$ then $f$ has a scrambled set consisting of transitive points $S=\bigcup_{n\geq 1}C_n$, where $C_n$ is a synchronously proximal Cantor set, i.e. $\liminf_{m\to\infty}\diam(f^m(C_n))=0$, and $\bigcup_{i=0}^{p-1}f^i(S)$ is dense in $X$.

In \cite{AGHSY2010}, another kind of chaos named uniform chaos was introduced and investigated by Akin at al. A subset $S$ is called a uniformly chaotic set if there are Cantor sets $C_1\subset C_2\cdots$ satisfying
\begin{enumerate}
	\item $S=\bigcup_{n\geq 1}C_n$ is a  recurrent and proximal subset of $X$;
	\item for any $n\geq 1$, $C_n$ is uniformly recurrent and uniformly proximal.
\end{enumerate}
Then $f$ is called uniformly chaotic if $f$ has a uniformly chaotic subset. A particular case of the main result in \cite{AGHSY2010} formulated the relationship between uniform chaos and Li-Yorke chaos.
For more information concerning chaos in topological dynamics the reader can refer to the survey \cite{LY2016}, where the development of chaos theory in topological dynamics is discussed in detail.

When the phase space is a vector space and the map is a linear map, some new characterizations of chaos can be obtained. Let $X$ be a Banach space over the field $\mathbb{K}$ and let $T\colon X\to X$ be a continuous linear operator.
Following \cite{B1988}, a vector $x\in X$ is said to be \emph{irregular} for $T$ if
\[
	\liminf_{n\to\infty} \Vert T^n x\Vert =0\text{ and }
	\limsup_{n\to\infty} \Vert T^n x\Vert =\infty.
\]
In \cite{BBMP2011}, Berm\'udez et al. characterized Li-Yorke chaos in terms of the existence of irregular vectors or the existence of a Li-Yorke scrambled pair for continuous linear operators on a Banach space.

The reader can also refer to \cite{BBMP2015} and \cite{JL2022} for more recent development of the study of Li-Yorke chaos for linear operators, where some new characterizations of Li-Yorke chaos and other kinds of chaos such as mean Li-Yorke chaos and distributional chaos are established.
In \cite{S2008}, Smith constructed an operator $T\colon \ell^2\to \ell^2$ such that  every non-zero vectors is irregular in the uniform sense, that is, there exist two sequences $\{p_n\}$ and $\{q_n\}$ in $\bbn$ such that
for any non-zero vector $x\in\ell^2$,
\[
	\lim_{n\to\infty} \Vert T^{p_n}x \Vert =0\text{ and }
	\lim_{n\to\infty} \Vert T^{q_n}x \Vert =\infty.
\]
Motivated by this example partially, we introduce the concept of uniform Li-Yorke chaos for linear operators.
To be specific, we say that a subset $S$ of $X$ is \emph{uniformly Li-Yorke scrambled} for $T$ if
there exist two sequences $\{p_n\}$ and $\{q_n\}$ in $\bbn$ such that
for any pair $x,y\in S$ with $x\neq y$, we have
\[
	\lim_{n\to\infty} \Vert T^{p_n}x -T^{p_n}y  \Vert =0\text{ and }
	\lim_{n\to\infty} \Vert T^{q_n}x -T^{q_n}y\Vert =\infty.
\]
Clearly, if $x$ is irregular then $\{\lambda x\colon \lambda\in \mathbb{K}\}$ is uniform Li-Yorke scrambled.
We say that an operator $T$ is \emph{densely uniformly Li-Yorke chaotic} if there exists a dense, uncountable, uniformly Li-Yorke scrambled set in $X$.

This paper focuses on the dense uniform Li-Yorke chaos for linear operators on a Banach space, and it is organized as follows. First of all, in Section \ref{sec:DULYmetric}, the problem is discussed in a  general setting, the dense uniform Li-Yorke chaos for continuous maps on a metric space is investigated, especially, some equivalent descriptions are established for dense uniform Li-Yorke chaos. In Section \ref{sec:DULYBanach}, we focus on linear systems defined on a Banach space. Due to the good properties of the Banach space, some sufficient conditions are formulated under which the operator is densely uniformly Li-Yorke chaotic.  We show that  there are plenty of densely uniformly Li-Yorke chaotic operators.
For unilateral backward weighted shifts and  bilateral backward weighted shifts on $\ell^p$, we show that Li-Yorke chaos is equivalent to dense uniform Li-Yorke chaos.

\section{Dense uniform Li-Yorke chaos for continuous maps on a metric space}\label{sec:DULYmetric}
Let $(X,d)$ be a metric space and $f\colon X\to X$ be a continuous map.
We say that a pair $(x,y)\in X\times X$ is \emph{Li-Yorke scrambled} if
\[
	\liminf_{n\to\infty} d(f^{n}(x) ,f^{n}(y))  =0
	\text{ and }
	\limsup_{n\to\infty} d(f^{n}(x), f^{n}(y))>0.
\]
A subset $S$ of $X$ with at least two points is called \emph{Li-Yorke scrambled} if every pair of two distinct points in $S$ forms a Li-Yorke scrambled pair.
We say that $f$ is \emph{(densely) Li-Yorke chaotic} if there exists a (dense) uncountable Li-Yorke scrambled subset of $X$.

Now we propose the concept of uniform Li-Yorke chaos, which is the starting point of this paper.
A subset $S$ of $X$ with at least two points is called \emph{uniformly Li-Yorke scrambled} for $f$ if
there exist two sequences $\{p_n\}$ and $\{q_n\}$ in $\bbn$ such that
for any pair $x,y\in S$ with $x\neq y$, we have
\[
	\lim_{n\to\infty} d(f^{p_n}(x) ,f^{p_n}(y))  =0\text{ and }
	\lim_{n\to\infty} d(f^{q_n}(x),  f^{q_n}(y)) =\infty.
\]
We say that $f$ is \emph{(densely) uniformly Li-Yorke chaotic} if there exists a (dense) uncountable uniformly Li-Yorke scrambled subset of $X$.

\begin{rem}
For a dynamical system $(X,T)$, if the metric $d$ on $X$ is bounded, then the system can not be  uniformly Li-Yorke chaotic.
However, there exists some Li-Yorke chaotic dynamical system such that any three points cannot be separated along a common sequence.

	By \cite{L2011}*{Theorem 4.21},
	for every continuous map $f\colon [0,1]\to[0,1]$,
	if the topological entropy of $f$ is zero then for every three points $x_1,x_2,x_3\in [0,1]$ one has
	\[
		\limsup_{n\to\infty} \min_{1\leq i<j\leq 3} d(f^n(x_i),f^n(x_j)) =0.
	\]
	There exists some continuous map $f\colon [0,1]\to[0,1]$ with zero topological entropy such that $f$ is Li-Yorke chaotic. By the above formula it is easy to see that this map $f$ is not uniformly Li-Yorke chaotic. 
\end{rem}

\begin{rem}
The proof of Theorem 4.2 in \cite{M1999} shows that there exists a $C^\infty$ map $f\colon \mathbb{R}\to\mathbb{R}$ such that the whole space $\mathbb{R}$ is uniformly Li-Yorke scrambled.
\end{rem}

Recall that a subset $A$ of $X$ is called a \emph{$G_\delta$ set} if it is the countable intersection of open sets,
a \emph{Cantor set}
if it is homeomorphic to the Cantor ternary set, and
a \emph{$\sigma$-Cantor set} if it is the countable union of Cantor sets.

For every $k\geq 2$, define
\[
	\mathrm{Prox}_k(f)=\biggl\{(x_1,x_2,\dotsc,x_k)\in X^k\colon
	\liminf_{n\to\infty} \max_{1\leq i<j\leq n} d(f^n(x_i),f^n(x_j)) =0\biggr\}
\]
and
\[
	D_k(f)=\biggl\{(x_1,x_2,\dotsc,x_k)\in X^k\colon
	\limsup_{n\to\infty} \min_{1\leq i<j\leq n} d(f^n(x_i),f^n(x_j)) =\infty\biggr\}.
\]

We have the following criterion for dense uniform Li-Yorke chaos.
\begin{thm} \label{thm:chaos-equi}
	Let $(X,d)$ be a complete separable metric space without isolated points and $f\colon X\to X$ be a continuous map. Then the following assertions are equivalent:
	\begin{enumerate}
		\item there exists a dense, uniformly Li-Yorke scrambled set in $X$;
		\item $f$ is densely uniformly Li-Yorke chaotic;
		\item there exists a dense, $\sigma$-Cantor, uniformly Li-Yorke scrambled set in $X$;
		\item for every $k\geq 2$, $\mathrm{Prox}_k(f)$ and $D_k(f)$ are dense in $X^k$.
	\end{enumerate}
\end{thm}
\begin{proof}
	By the definitions of dense uniform Li-Yorke chaos and uniformly Li-Yorke scrambled set and noticing the constructions of $\mathrm{Prox}_k(f)$ and $D_k(f)$ we know that (3)$\Rightarrow$(2)$\Rightarrow$(1)$\Rightarrow$(4) are clear.

	We only need to prove (4)$\Rightarrow$(3).
	As $X$ is separable, pick a countable topological basis $\{O_n\}_{n=1}^\infty$ of $X$.
	Let $a_0=0$ and $a_{n+1}=2 a_n+2$ for $n\geq 1$.
	We are going to construct two sequences $\{p_n\}$ and $\{q_n\}$ in $\bbn$, and nonempty open subsets $U^{(n)}_i$ in $X$ with $n=1,2,\dotsc$ and $1\leq i\leq a_n$ such that
	\begin{enumerate}[(i)]
		\item $\overline{U^{(n+1)}_i}\cup \overline{ U^{(n+1)}_{a_n+i}}\subset U^{(n)}_i$ for $1\leq i\leq a_n$, and $\overline{U^{(n+1)}_{2a_n+1}}\cup \overline{U^{(n+1)}_{2a_n+2}}\subset O_{n+1}$;
		\item for every $n\in\bbn$, $\Bigl\{\overline{U^{(n)}_i}\Bigr\}_{1\leq i\leq a_n}$ are pairwise disjoint;
		\item $\diam(U^{(n)}_i)<\frac{1}{n}$ for $n=1,2,\dotsc$, and $1\leq i\leq a_n$;
		\item for every $x_i\in U^{(n)}_i$ with $1\leq i\leq a_n$,
		      \[
			      \max_{1\leq i<j\leq a_n} d(f^{p_n}(x_i),f^{p_n}(x_j))<\tfrac{1}{n};
		      \]
		\item for every $x_i\in U^{(n)}_i$ with $1\leq i\leq a_n$,
		      \[
			      \min_{1\leq i<j\leq a_n} d(f^{q_n}(x_i),f^{q_n}(x_j))>n.
		      \]
	\end{enumerate}
	For $n=1$, as $X$ has no isolated points,
	pick two nonempty open subsets $V^{(1)}_1$ and $V^{(1)}_2$ of $O_1$ such that
	\begin{enumerate}[(i)]
		\item $\overline{V^{(1)}_1}\cup \overline{V^{(1)}_2}\subset O_1$;
		\item $\overline{V^{(1)}_1}\cap \overline{V^{(1)}_2}=\emptyset$;
		\item $\diam(V^{(1)}_i)<1$ for $i=1,2$.
	\end{enumerate}
	As $\mathrm{Prox}_{a_1}(f)$ is dense in $X^{a_1}$,
	pick $y_i\in V^{(1)}_i$ for $i=1,2$ and $p_1\in\bbn$
	such that
	\[
		d(f^{p_1}(y_1),f^{p_1}(y_2))<1.
	\]
	By the continuity of $f^{p_1}$, there exist
	nonempty open subsets $W^{(1)}_i$ of $V^{(1)}_i$ for $i=1,2$ such that
	\begin{enumerate}
		\item[(iv)]for every $x_i\in W^{(1)}_i$ with $i=1,2$,
		\[
			d(f^{p_1}(x_1),f^{p_1}(x_2))<1.
		\]
	\end{enumerate}
	As $D_{a_1}(f)$ is dense in $X^{a_1}$,
	pick $z_i\in W^{(1)}_i$ for $i=1,2$ and $q_1\in\bbn$
	such that
	\[
		d(f^{q_1}(z_1),f^{q_1}(z_2))>1.
	\]
	By the continuity of $f^{q_1}$, there exist
	nonempty open subsets $U^{(1)}_i$ of $W^{(1)}_i$ for $i=1,2$ such that
	\begin{enumerate}
		\item[(v)]
			for every $x_i\in U^{(1)}_i$, $i=1,2$,
			\[
				d(f^{q_1}(x_1),f^{q_1}(x_2))>1.
			\]
	\end{enumerate}

	Assume that $\{a_j\}_{j=1}^n$, $\{p_j\}_{j=1}^n$ and $\{q_j\}_{j=1}^n$ in $\bbn$, and nonempty open subsets $U^{(j)}_i$ in $X$ with $j=1,2,\dotsc,n$ and $1\leq i\leq a_n$ have been constructed satisfying (i)--(v).
	Let $a_{n+1}=2a_n+2$.
	As $X$ has no isolated points,
	pick  nonempty open subsets $V^{(n+1)}_i$ for $i=1,2,\dotsc,a_{n+1}$ such that
	\begin{enumerate}[(i)]
		\item $\overline{V^{(n+1)}_i}\cup\overline{ V^{(n+1)}_{a_n+i}}\subset U^{(n)}_i$ for $1\leq i\leq a_n$, and $\overline{V^{(n+1)}_{2a_n+1}}\cup\overline{V^{(n+1)}_{2a_n+2}}\subset O_{n+1}$;
		\item $\Bigl\{\overline{V^{(n+1)}_i}\Bigr\}_{1\leq i\leq a_{n+1}}$ are pairwise disjoint;
		\item $\diam(V^{(n+1)}_i)<\frac{1}{n+1}$ for   $1\leq i\leq a_{n+1}$.
	\end{enumerate}
	As $Prox_{a_{n+1}}(f)$ is dense in $X^{a_{n+1}}$,
	pick $y_i\in  V^{(n+1)}_i$ for $i=1,2,\dotsc,a_{n+1}$ and $p_{n+1}\in\bbn$ with $p_{n+1}>p_{n}$
	such that
	\[
		\max_{1\leq i<j\leq a_{n+1}} d(f^{p_{n+1}}(y_i),f^{p_{n+1}}(y_j))<\tfrac{1}{n+1}.
	\]
	By the continuity of $f^{p_{n+1}}$, there exist
	nonempty open subsets $W^{(n+1)}_i$ of $V^{(n+1)}_i$ for $i=1,2,\dotsc,a_{n+1}$ such that
	\begin{enumerate}
		\item[(iv)] for every $x_i\in W^{(n+1)}_i$ with $i=1,2,\dotsc,a_{n+1}$,
			\[
				\max_{1\leq i<j\leq a_{n+1}} d(f^{p_{n+1}}(x_i),f^{p_{n+1}}(x_j))<\tfrac{1}{n+1}.
			\]
	\end{enumerate}

	As $D_{a_{n+1}}(f)$ is dense in $X^{a_{n+1}}$,
	pick $z_i\in W^{(n+1)}_i$ for $i=1,2,\dotsc,a_{n+1}$ and $q_{n+1}\in\bbn$
	with $q_{n+1}>q_n$
	such that
	\[
		\min_{1\leq i<j\leq a_{n+1}} d(f^{q_{n+1}}(z_i),f^{q_{n+1}}(z_j))>n+1.
	\]
	By the continuity of $f^{q_{n+1}}$, there exist
	nonempty open subsets $U^{(n+1)}_i$ of $W^{(n+1)}_i$ for $i=1,2,\dotsc,a_{n+1}$ such that
	\begin{enumerate}
		\item[(v)]
			for every $x_i\in U^{(n+1)}_i$ with  $i=1,2,\dotsc,a_{n+1}$,
			\[
				\min_{1\leq i<j\leq a_{n+1}} d(f^{q_{n+1}}(x_i),f^{q_{n+1}}(x_j))>n+1.
			\]
	\end{enumerate}
	By the induction, we have construed the desired sequences.
	For every $k\in\bbn$, let
	\[
		S_k = \bigcap_{n=k}^\infty \bigcup_{i=1}^{2^{n-k}a_k} \overline{U^{(n)}_i}.
	\]
	Then $S_1\subset S_2\subset \dotsb$ and each $S_k$ is a Cantor set.
	Finally let $S=\bigcup_{k =1}^\infty S_k$.
	Then it is easy to see that $S$ is a dense, $\sigma$-Cantor, uniformly Li-Yorke scrambled set in $X$.
\end{proof}

By the proof of Theorem~\ref{thm:chaos-equi}, we have the following two consequences.

\begin{coro}\label{coro:prox-k-dense}
	Let $(X,d)$ be a complete separable metric space without isolated points and $f\colon X\to X$ be a continuous map.
	Then the following assertions are equivalent:
	\begin{enumerate}
		\item for each $k\geq 2$, $\mathrm{Prox}_k(f)$ is dense in $X^k$;
		\item there exists a sequence $\{p_n\}$  in $\bbn$ and a dense $\sigma$-Cantor subset $S$ of $X$ such that
		      for any pair $x,y\in S$
		      \[
			      \lim_{n\to\infty} d(f^{p_n}(x) ,f^{p_n}(y))  =0.
		      \]
	\end{enumerate}
\end{coro}

\begin{coro}\label{coro:Dk}
	Let $(X,d)$ be a complete separable metric space without isolated points and $f\colon X\to X$ be a continuous map.
	Then the following assertions are equivalent:
	\begin{enumerate}
		\item for each $k\geq 2$, $D_k(f)$ is dense in $X^k$;
		\item there exists a sequence $\{q_n\}$  in $\bbn$ and a dense $\sigma$-Cantor subset $S$ of $X$ such that
		      for any pair $x,y\in S$ with $x\neq y$
		      \[
			      \lim_{n\to\infty} d(f^{q_n}(x) ,f^{q_n}(y)) =\infty.
		      \]
	\end{enumerate}
\end{coro}

Let $(X,d)$ be a complete separable metric space and $f\colon X\to X$ be a continuous map.
We say that $f$ is \emph{transitive} if for every two nonempty open subsets $U$ and $V$ of $X$ there exists $n\in\bbn$ such that $U\cap f^{-n}(V)\neq\emptyset$, \emph{weakly mixing} if the product map $f\times f\colon X\times X\to X\times X$, $(x,y)\mapsto (f(x),f(y))$, is transitive,
and \emph{strongly mixing} if for every two nonempty open subsets $U$ and $V$ of $X$ there exists $N\in\bbn$ such that $U\cap f^{-n}(V)\neq\emptyset$ for all $n\geq N$.
The well-known Furstenberg intersection lemma states that if $f\colon X\to X$ is weakly mixing then for any $k\geq 2$ the $k$th  product system $(X^k,f\times f\times \dotsb \times f)$ is transitive (see e.g. \cite{GP2011}*{Theorem 1.51}).

For a point $x\in X$, the \emph{orbit} of $x$ is the set
\[
	\mathrm{Orb}_f(x):=\{f^n(x)\colon n\geq 0\}.
\]
If the orbit of $x$ is dense in $X$, then we say that $x$ is a \emph{transitive point}. By the Birkhoff transitive theorem, if $(X,d)$ is a complete separable metric space without isolated points, then a map $f\colon X\to X$ is transitive if and only if the collection of transitive points of $f$ is a dense $G_\delta$ subset of $X$ (see e.g. \cite{GP2011}*{Theorem 1.16}).

\begin{rem}
	Let $(X,d)$ be a complete separable metric space without isolated points and $f\colon X\to X$ be a continuous map.
	In \cite{HY2002} Huang and Ye proved that if $f$ is transitive and has a fixed point $z\in X$ then it is densely Li-Yorke chaotic.
	Moreover, \cite{M2004} Mai showed that if $f$ is transitive and has a fixed point $z\in X$ then the conditions in Corollary \ref{coro:prox-k-dense} are satisfied.
	Let us explain this for completeness.
	Fix a transitive point $x\in X$. Then there exists a sequence $\{p_n\}$ in $\bbn$ such that $\lim_{n\to\infty} f^{p_n}(x)=z$.
	By the continuity of $f$, for every $i\in\bbn$, $\lim_{n\to\infty} f^{p_n}(f^i(x))=f^{i}(z)=z$.
	As $x$ is a transitive point,
	for every $k\geq 2$,
	the set $\{(y_1,\dotsc,y_k)\colon y_i\in \mathrm{Orb}_f(x)\}$ is dense in $X^k$ and contained in $\mathrm{Prox}_k(f)$.
\end{rem}

It is clear that if $f\colon X\to X$ is  uniformly Li-Yorke chaotic then the diameter of $X$ must be infinite.
The following result reveals that under this condition weak mixing implies dense uniform Li-Yorke chaos.

\begin{prop}\label{prop:wm-DULYC}
	Let $(X,d)$ be a complete separable metric space and $f\colon X\to X$ be a continuous map.
	If $\diam(X)=\infty$ and $f$ is weakly mixing then $f$ is densely uniformly Li-Yorke chaotic.
\end{prop}
\begin{proof}
	As the diameter of $X$ is infinite, $X$ is an infinite set.
	Then $X$ has no isolated points, since $f$ is weakly mixing.
	For each $k\in\bbn$, let $\mathrm{Trans}_k(f)$ be the collection of transitive points for the product system $(X^k,f\times f\times \dotsb \times f)$. Then $\mathrm{Trans}_k(f)$ is a dense $G_\delta$ subset of $X^k$.
	By Theorem~\ref{thm:chaos-equi},  it is sufficient to show that $\mathrm{Trans}_k(f)\subset \mathrm{Prox}_k(f)\cap D_k(f)$.

	Fix a point $z\in X$. For every $(x_1,\dotsc,x_k)\in \mathrm{Trans}_k(f)$, there exists a sequence $\{p_n\}$ such that
	$f^{p_n}(x_j)\to z$ as $n\to\infty$ for $j=1,\dotsc,k$.
	Then
	\[
		\liminf_{n\to\infty} \max_{1\leq i<j\leq n} d(f^n(x_i),f^n(x_j)) =0,
	\]
	that is $(x_1,\dotsc,x_k)\in \mathrm{Prox}_k(f)$.

	Fix $M>0$.
	As $\diam(X)=\infty$, there exist $z_1,z_2\in X$ such that
	$d(z_1,z_2)\geq M$.
	By $\diam(X)=\infty$ again, two open balls $B(z_1,M)$ and
	$ B(z_2,M)$ do not cover $X$. Then there exists a point $z_3\in X$ such that $d(z_3,z_1)\geq M$ and $d(z_3,z_2)\geq M$.
	By recursive construction, we get a point $(z_1,\dotsc,z_k)\in X^k$ such that $d(z_i,z_j)\geq M$ for $1\leq i<j\leq k$.
	For every $(x_1,\dotsc,x_k)\in \mathrm{Trans}_k(f)$, there exists a sequence $\{q_n\}$ such that
	$f^{q_n}(x_j)\to z_j$ as $n\to\infty$ for $j=1,\dotsc,k$.
	Then
	\[
		\liminf_{n\to\infty} \max_{1\leq i<j\leq n} d(f^n(x_i),f^n(x_j)) \geq M.
	\]
	By the arbitrariness of $M$, we get \[
		\liminf_{n\to\infty} \max_{1\leq i<j\leq n} d(f^n(x_i),f^n(x_j)) =\infty,
	\]
	that is $(x_1,\dotsc,x_k)\in D_k(f)$.
\end{proof}

Let $(X,d)$ be a complete separable metric space.
A subset $A$ of $X$ is called \emph{residual} if it contains a dense $G_\delta$ of $X$, \emph{of first category} if it is the complementary set of a residual set,
and \emph{of second category} if it is not of first category.
It is clear that a subset of $X$ is of second category if and only if it intersects every residual subset of $X$.
The following result shows that a strongly mixing map only can have ``small'' uniformly Li-Yorke scrambled sets.

\begin{prop}
	Let $(X,d)$ be a complete separable metric space and $f\colon X\to X$ be a continuous map.
	If  $f$ is strongly mixing, then every uniformly Li-Yorke scrambled set of $f$ must be of first category.
\end{prop}

\begin{proof}
	We prove this result by contradiction.
	Assume that $f$ has a   uniformly Li-Yorke scrambled subset $S$ of second category.
	Then there exists a sequence $\{p_n\}$ in $\bbn$ such that
	for any pair $x,y\in S$ we have
	\[
		\lim_{n\to\infty} d(f^{p_n}(x) ,f^{p_n}(y))  =0.
	\]
	It is clear that the product map $f\times f\colon X\times X\to X\times X$ is also strongly mixing.
	Pick two nonempty open subsets $U$ and $V$ of $X$ with $U\cap V=\emptyset$.
	Let
	\[
		A=\bigcap_{N=1}^\infty \bigcup_{n=N}^\infty (f\times f)^{-p_n} (U\times V).
	\]
	As $f\times f$ is strongly mixing, it is easy to see that $A$ is a dense $G_\delta$ subset of $X\times X$.
	By the Claim in the proof of \cite{K1991}*{Theorem~2},
	$S\times S$ is of second category in $X\times X$. Then $A \cap (S\times S)\neq\emptyset$.
	Pick a point $(x,y)\in A\cap (S\times S)$. There is an infinite subsequence $\{p_{n_j}\}$ of $\{p_n\}$ such that $f^{p_{n_j}}(x)\in U$ and $f^{p_{n_j}}(y)\in V$ for all $j\in\bbn$. This contradicts to the fact $\lim_{n\to\infty} d(f^{p_n}(x) ,f^{p_n}(y)) = 0$.
\end{proof}

\begin{rem}
	Let $(X,d)$ be a complete separable metric space, $f\colon X\to X$ be a continuous map and $k\in\bbn$.
	If a subset $S$ of $X$ is uniformly Li-Yorke scrambled for $f^k$, then it is clear that it is also true for $f$.
	If in addition $f$ is Lipschitz continuous then it is not hard to check that the inverse is true.
\end{rem}

\section{Dense uniform Li-Yorke chaos for operators on a Banach space}\label{sec:DULYBanach}

In this section unless otherwise specified, let $X$ be a separable Banach space over the field $\mathbb{K}$ ($\mathbb{R}$ or $\mathbb{C}$)
with norm $\Vert\cdot\Vert$ and $T\colon X\to X$ be a bounded linear operator from $X$ to itself. A transitive point in $X$ is usually called a \emph{hypercyclic vector} for $T$, and we say that $T$ is \emph{hypercyclic} if it has some hypercyclic vector.

Let $\mathcal{B}(X)$ be the collection of bounded operators from $X$ to itself.
The \emph{strong operator topology} (SOT for short) on $\mathcal{B}(X)$
is defined as follows: any $T\in \mathcal{B}(X)$ has a neighborhood basis consisting of  sets of the form
\[
	V_{x_1,\dotsc,x_n,\varepsilon,T}=\{
	S\in \mathcal{B}(X)\colon \|Sx_i-Tx_i\|<\varepsilon,\ i=1,\dotsc,n\},
\]
where $x_1,\dotsc,x_n\in X$ and $\varepsilon>0$.
Denote by $\mathrm{DULYC}(X)$ the collection of all densely uniformly Li-Yorke chaotic operators in $\mathcal{B}(X)$.
Combining Theorem 2 and Remark ii) in page 204 of \cite{BC2003}  the collection of weakly mixing operators on $X$ is SOT-dense in $\mathcal{B}(X)$.
By Proposition~\ref{prop:wm-DULYC}, every weakly mixing operator is densely uniformly Li-Yorke chaotic. So we have the following consequence.

\begin{prop}\label{prop:SOT-dense}
	Let $X$ be an infinite-dimensional separable Banach space. Then $\mathrm{DULYC}(X)$ is  SOT-dense  in  $\mathcal{B}(X)$.
\end{prop}

If in addition $X$ is a separable Hilbert space.
The \emph{strong$^*$ operator topology} (SOT$^*$ for short) on $\mathcal{B}(X)$
is defined as follows: any $T\in \mathcal{B}(X)$ has a neighborhood basis consisting of  sets of the form
\begin{align*}
	V^*_{x_1,\dotsc,x_n,\varepsilon,T}=\{
	S\in \mathcal{B}(X)\colon \|Sx_i-Tx_i\|<\varepsilon, \|S^*x_i-T^*x_i\|<\varepsilon,\ i=1,\dotsc,n\},
\end{align*}
where $x_1,\dotsc,x_n\in X$, $\varepsilon>0$ and $T^*$ is the adjoint operator of $T$.
It is clear that the  strong operator topology is coarser than the strong$^*$ operator topology.
For $M>0$, the closed ball $\mathcal{B}_M(X)=\{T\in \mathcal{B}(X)\colon \|T\|\leq M\}$ with respect to strong operator topology or strong$^*$ operator topology becomes a Polish (separable and completely metrizable) space, see e.g. \cite{K1995}*{Exercise I.3.4(5)}.

\begin{prop}
	Assume that $X$ is an infinite-dimensional separable Hilbert space and $M>1$.
	Then $\mathrm{DULYC}(X)\cap \mathcal{B}_M(X)$ is a SOT$^*$-dense $G_\delta$ subset of $\mathcal{B}_M(X)$.
\end{prop}
\begin{proof}
We first show the following claim.

\medskip 
\noindent\textbf{Claim: }
$\mathrm{DULYC}(X)\cap \mathcal{B}_M(X)$ is a $G_\delta$ subset of $\mathcal{B}_M(X)$ with respect to the strong operator topology.

\begin{proof}[Proof of The Claim]
Let $\mathcal{U}$ be a countable topological basis of $X$. Without loss of generality, assume that open sets in $\mathcal{U}$ are nonempty.
	Fix $k\geq 2$, $N\in\bbn$ and $U_1,\dotsc,U_k\in \mathcal{U}$.
	Let
	\begin{align*}
		P(U_1,\dotsc,U_k,N)=\{ T\in \mathcal{B}_M(X)\colon
		 & \exists x_i\in U_i \text{ for } i=1,\dotsc,k \text{ and } n\geq N \\
		 & \text{ s.t. }
		\Vert T^n(x_i-x_j)\Vert <\tfrac{1}{N}  \text{ for } 1\leq i< j\leq k\}
	\end{align*}
	and
	\begin{align*}
		D(U_1,\dotsc,U_k,N)=\{ T\in \mathcal{B}_M(X)\colon
		 & \exists x_i\in U_i \text{ for } i=1,\dotsc,k \text{ and } n\geq N \\
		 & \text{ s.t. }
		\Vert T^n(x_i-x_j)\Vert >N \text{ for } 1\leq i<j\leq k\}.
	\end{align*}
	For every $n\in\bbn$, the map $(\mathcal{B}_M(X),\mathrm{SOT})\to (\mathcal{B}(X),\mathrm{SOT})$, $T\mapsto T^n$ is continuous, see e.g. \cite{GMM2021}*{Lemma 2.1}. It is easy to check that
	$P(U_1,\dotsc,U_k,N)$ and $D(U_1,\dotsc,U_k,N)$ are open in $(\mathcal{B}_M(X), \mathrm{SOT})$.
	Now by Theorem~\ref{thm:chaos-equi},
	\[
		\mathrm{DULYC}(X)\cap \mathcal{B}_M(X)=\bigcap_{k=2}^\infty \bigcap_{U_1,\dotsc,U_k\in \mathcal{U}} \bigcap_{N=1}^\infty P(U_1,\dotsc,U_k,N)  \cap D(U_1,\dotsc,U_k,N).
	\]
	This shows that $\mathrm{DULYC}(X)\cap \mathcal{B}_M(X)$ is a $G_\delta$ subset of $\mathcal{B}(X)$.
 \end{proof}
 
Since the strong$^*$ operator topology is finer than the strong operator topology, by the Claim 
$\mathrm{DULYC}(X)\cap \mathcal{B}_M(X)$ is a  $G_\delta$ subset of $\mathcal{B}_M(X)$ with respect to the strong$^*$ operator topology.
By \cite{GMM2021}*{Proposition 2.16} the collection of weak mixing operators in $\mathcal{B}_M(X)$ is SOT$^*$-dense in $\mathcal{B}_M(X)$.
	Hence, $\mathrm{DULYC}(X)\cap \mathcal{B}_M(X)$ is a SOT$^*$-dense $G_\delta$ subset of $\mathcal{B}_M(X)$.
\end{proof}

Compared with Corollary~\ref{coro:prox-k-dense}, by the linearity of $T$ we have the following simple characterization.
\begin{lem}\label{lem:T-prox-dense-equi}
	Let $T\in\mathcal{B}(X)$.
	Then for every $k\geq 2$, $\mathrm{Prox}_k(f)$ is  dense in $X^k$
	if and only if there exist a sequences $\{p_n\}$ in $\bbn$
	and a dense subset $K$ of $X$ such that  for every $x\in K$,
	$\lim_{n\to\infty} T^{p_n}x =0 $.
\end{lem}
\begin{proof}
	Suppose that $\mathrm{Prox}_k(f)$ is dense in $X^k$ for every $k\geq 2$. By Theorem \ref{thm:chaos-equi} we know that there exist a sequence $\{p_n\}$ in $\bbn$ and a dense subset $S$ of $X$ such that for any pair $x, y\in S$ we have
	\[
		\lim_{n\to\infty}\Vert T^{p_n}x-T^{p_n}y\Vert=0.
	\]
	Fix $y\in S$ and let $K=S-y$.
	Then $K$ is a dense subset of $X$ and for any $x\in K$,
	\[
		\lim_{n\to\infty} T^{p_n}x=0.
	\]

	Now suppose that there exist a sequence $\{p_n\}$ in $\bbn$ and a dense subset $K$ of $X$ such that for every $x\in K$, $\lim_{n\to\infty}T^{p_n}x= 0$.
	Then for any $x_1,x_2\in K$ with $x_1\not= x_2$ we have
	\[
		\lim_{n\to\infty}\Vert T^{p_n}x_1-T^{p_n}x_2\Vert
		\leq\lim_{n\to\infty}\Vert T^{p_n}x_1\Vert+\lim_{n\to\infty}\Vert T^{p_n}x_2\Vert=0,
	\]
	which completes the proof the lemma.
\end{proof}

Let $(X,d)$ be a metric space and $f\colon X\to X$ be a continuous map.
We say that $f$ has \emph{sensitive dependence on initial conditions} or is \emph{sensitive} briefly if there exists some $\delta>0$ such that for every $x\in X$ and $\varepsilon>0$, there exists some $y\in X$ with $d(x,y)<\varepsilon$ and some $n\in\bbn$ such that
$d(f^n(x),f^n(y))>\delta$.

Using the Banach-Steinhaus theorem, we have the following characterization of sensitive operators on a Banach space,
see e.g. \cite{GP2011}*{Exercise 2.3.1} or
\cite{JL2022}*{Proposition 3.29}.

\begin{lem}\label{lem:T-sensitive-eq}
	Let $T\in\mathcal{B}(X)$.
	Then the following assertions are equivalent:
	\begin{enumerate}
		\item $T$ is sensitive;
		\item $\sup_n\Vert T^n\Vert=\infty$;
		\item there exists a vector $x\in X$ such that $\limsup_{n\to\infty} \Vert T^n x\Vert =\infty$;
        \item the set of all vectors $x\in X$ satisfying 
        $\limsup_{n\to\infty} \Vert T^n x\Vert =\infty$ is residual in $X$.
	\end{enumerate}
\end{lem}

Let $T\in\mathcal{B}(X)$.
A vector subspace $Y$ of $X$ is called an \emph{irregular manifold} for $T$ if
every non-zero vector  $y\in Y$ is irregular for $T$.
According to \cite{BBMP2015}*{Theorem 31}, we have the following sufficient criterion for the existence
of a dense irregular manifold.  
\begin{thm}\label{thm:dense-irregular-mainfold}
Let $T\in\mathcal{B}(X)$. 
If $T$ is sensitive and there exists a sequence $\{p_n\}$ in $\bbn$ and a dense subset $K$ of $X$ such that  for every $x\in K$,
$\lim_{n\to\infty} T^{p_n}x =0 $,
then $T$ admits a dense irregular manifold.
\end{thm}

Combing Lemmas~\ref{lem:T-prox-dense-equi} and~\ref{lem:T-sensitive-eq}, Theorems~\ref{thm:chaos-equi} and~\ref{thm:dense-irregular-mainfold},
we have the following consequence.

\begin{coro}\label{cor:DULYC-dense-irr-manifold}
Let $T\in\mathcal{B}(X)$. 
If $T$ is densely uniformly Li-Yorke chaotic, then $T$  admits a dense irregular manifold.
\end{coro}

We have the following sufficient conditions of dense uniform Li-Yorke chaos for operators.

\begin{thm}\label{thm:T-asy-sens}
	Let $T\in\mathcal{B}(X)$.
	If $T$ is sensitive and $\{x\in X\colon \lim_{n\to\infty} T^nx =0\}$ is dense in $X$,
	then it is densely uniformly Li-Yorke chaotic.
\end{thm}
\begin{proof}
	By Lemma~\ref{lem:T-prox-dense-equi}, we have for every $k\geq 2$, $\mathrm{Prox}_k(f)$ is  dense in $X^k$.
	By Theorem \ref{thm:chaos-equi}, it is sufficient to show that for every $k\geq 2$, $D_k(f)$ is  dense in $X^k$.
	Fix $k\geq 2$ and nonempty open subsets $U_1,\dotsc,U_k$ of $X$.
	Pick points $x_i\in U_i$ such that $\lim_{n\to\infty} T^n x_i=0$ for $i=1,\dotsc,k$.
	By Lemma~\ref{lem:T-sensitive-eq}, pick a point $z\in X$ such that
	$\limsup_{n\to\infty} \Vert T^n z\Vert =\infty$.
	Pick distinct scalars $\alpha_1,\dotsc,\alpha_k\in\mathbb{K}$ such that $x_i + \alpha_i z\in U_i$ for $i=1,\dotsc,k$. Then $(x_1+\alpha_1z,\dotsc,x_k+\alpha_kz)\in D_k(f)\cap U_1\times \dotsc\times U_k$.
	This shows that $D_k(f)$ is  dense in $X^k$.
\end{proof}

\begin{rem}
	Let $T\in\mathcal{B}(X)$.
	Recall that the \emph{kernel} of $T$ is $\mathrm{ker}(T)=\{x\in X\colon Tx=0\}$ and \emph{generalized  kernel of $T$} is $\bigcup_{n\in\bbn} \mathrm{ker}(T^n)$.
	According to Theorem~\ref{thm:T-asy-sens}, if the generalized kernel of $T$ is dense in $X$ then $T$ is sensitive if and only if it is densely uniformly Li-Yorke chaotic.
\end{rem}

To apply Theorem \ref{thm:T-asy-sens}, we will consider the unilateral backward weighted shift $B_w\colon X\to X$, $(x_1,x_2,\dotsc)\mapsto (w_2x_2,w_3x_3,\dotsc)$ on $X=\ell^p(\bbn)$,
$1\leq p<\infty$ or $X=c_0(\bbn)$,
where $w=\{w_j\}_{j\in\bbn}$ is a bounded sequence of non-zero weights.
It is easy to verify that for each $n\in\bbn$, $\Vert B_w^n\Vert =\sup_{k\in\bbn} \prod_{j=k+1}^{n+k+1}|w_j|$.
It is shown in \cite{BBMP2011}*{Proposition 27} that
for a unilateral backward weighted shift $B_w$ it is Li-Yorke chaotic if and only if $\sup_{k,n\in\bbn}\prod_{j=k+1}^{n+k+1} |w_j|=\infty$.
The generalized  kernel of $B_w$ is the collection of sequences with only finite non-zero coordinates, which is dense in $X$.
By Theorem \ref{thm:T-asy-sens}, we have the following result.

\begin{coro}\label{coro:unilateral-bws}
	Let $B_w$ be a unilateral backward weighted shift on $X=\ell^p(\bbn)$, $1\leq p<\infty$ or $X=c_0(\bbn)$ with the weight sequence $w=\{w_j\}_{j\in\bbn}$.
	Then the following conditions are equivalent:
	\begin{enumerate}
		\item $\sup_{k,n\in\bbn}\prod_{j=k+1}^{n+k+1} |w_j|=\infty$;
		\item there exists a Li-Yorke scrambled pair;
		\item $B_w$ is Li-Yorke chaotic;
          \item $B_w$ admits a dense  irregular manifold;
		\item $B_w$ is densely uniformly Li-Yorke chaotic.
	\end{enumerate}
\end{coro}

Now we consider the bilateral backward weighted shift $B_w\colon X\to X$, $(\dotsc,x_{-1},\allowbreak [x_0],  x_1,\dotsc)\mapsto
	(\dotsc,w_0x_{0},[w_1x_1], w_2x_2,\dotsc)$ on $X=\ell^p(\bbz)$,
$1\leq p<\infty$ or $X=c_0(\bbz)$,
where $w=\{w_j\}_{j\in\bbz}$ is a bounded sequence of non-zero weights.
It is easy to verify that for each $n\in\bbn$, $\Vert B_w^n\Vert =\sup_{k\in\bbz} \prod_{j=k}^{n+k}|w_j|$.
Note that in this case the generalized kernel of $B_w$ is trivial, so we can not apply Theorem~\ref{thm:T-asy-sens} directly.
We have the following characterization of dense uniform Li-Yorke chaos for bilateral backward weighted shifts.

\begin{thm}\label{thm:bi-backward-LY-chaos}
	Let $B_w$ be a bilateral backward weighted shift on $X=\ell^p(\bbz)$, $1\leq p<\infty$ or $X=c_0(\bbz)$ with the weight sequence $w=\{w_j\}_{j\in\bbn}$.
	Then the following conditions are equivalent:
	\begin{enumerate}
		\item $\liminf\limits_{n\to\infty} \prod_{j=-n+1}^0 |w_j|=0$ and $\sup_{k\in\bbz,n\in\bbn}\prod_{j=k}^{k+n} |w_j|=\infty$;
		\item there exists a Li-Yorke scrambled pair;
		\item $B_w$ is Li-Yorke chaotic;
		\item $B_w$ admits a dense irregular manifold;
		\item $B_w$ is densely uniform Li-Yorke chaotic.
	\end{enumerate}
\end{thm}

\begin{proof}
	For simplicity, we only prove the case $X=\ell^1(\bbz)$, the proofs of other cases are similar.

	(5)$\Rightarrow$(3)$\Rightarrow$(2) are obvious 
 and (5)$\Rightarrow$(4) follows from Corollary~\ref{cor:DULYC-dense-irr-manifold}.

	(2)$\Rightarrow$(1). Let $(x,y)$ be a Li-Yorke scrambled pair and let $z=x-y$.
	Then $z\neq 0$. Pick $k\in\bbz$ such that $z_k\neq 0$.
	For every $n\in\bbn$,
	$\Vert B_w^n z\Vert \geq |z_k|\cdot \prod_{j=k-n+1}^k |w_j |$.
	As $\liminf_{n\to\infty} \Vert B_w^n z \Vert =0$,
	one has $\liminf_{n\to\infty} |\prod_{j=k-n+1}^k w_j | =0$.
	Since the weights are not equal to zero, this is equivalent to
	$\liminf\limits_{n\to\infty} \prod_{j=-n+1}^0 |w_j|=0$.

	To show that $\sup_{k\in\bbz,n\in\bbn}\prod_{j=k}^{k+n} |w_j|=\infty$, by Lemma~\ref{lem:T-sensitive-eq} it is sufficient to show that $B_w$ is sensitive.
	Let $\delta = \limsup_{n\to\infty} \Vert B_w^n z\Vert$. Then $\delta>0$.
	If $\delta=\infty$, then by  Lemma~\ref{lem:T-sensitive-eq} $B_w$ is sensitive. Now assume that $\delta<\infty$.
	Fix $x\in X$ and $\varepsilon>0$.
	There exist $n_1,n_2\in\bbn$ with $n_1>n_2$ such that
	$\Vert B_w^{n_1}z\Vert <\varepsilon$ and $\Vert B_w^{n_2}z\Vert >\frac{1}{2}\delta$.
	Then $\Vert x-(x+B_w^{n_1}z)\Vert<\varepsilon$ and
	$\Vert B_w^{n_2-n_1}x - B_w^{n_2-n_1}(x+B_w^{n_1}z)\Vert = \Vert B_w^{n_2}z\Vert >\frac{1}{2}\delta$. This implies that $B_w$ is sensitive.

	(1)$\Rightarrow$(5).
	Let $\{e_i\}_{i\in\bbz}$ be the canonical basis of $\ell^1(\bbz)$.
	Note that by the definition, one has $B_w e_i =w_ie_{i-1}$ for all $i\in\bbz$.
	There exists a sequence $\{p_n\}$ in $\bbn$ such that $\lim\limits_{n\to\infty} \prod_{j=-p_n+1}^0 |w_j|=0$.
	Then for every $s\geq 1$,
	\[
		\lim\limits_{n\to\infty}  \Vert B_w^{p_n+s} e_k\Vert =
		\lim\limits_{n\to\infty} \prod_{j=-p_n+1}^s |w_j|=0.
	\]
	If $x=\sum_{j=-s}^s\alpha_j e_j$, then
	\begin{align*}
		\lim_{n\to\infty}  \Vert B_w^{p_n+s} x\Vert
		 & \leq \lim_{n\to\infty} \sum_{j=-s}^s |\alpha_j| \cdot \Vert B_w^{p_n+s} e_j\Vert
		=
		\lim_{n\to\infty}\sum_{j=-s}^s |\alpha_j| \cdot \Vert B_w^{p_n+s} (B_w^{s-j} e_s)\Vert
		\\
		 & \leq \lim_{n\to\infty} \sum_{j=-s}^s |\alpha_j| \cdot
		\Vert B_w^{s-j}\Vert \cdot
		\Vert B_w^{p_n+s} e_s \Vert=0.
	\end{align*}
	Fix $k\geq 2$ and nonempty open subsets $U_1,\dotsc,U_k$ of $X$.
	There exists $s\in \bbn$ and $x_i\in U_i$ such that $x_i$ can be expressed as
	$\sum_{j=-s}^s\alpha_j^{(i)} e_j$ for $i=1,\dotsc,k$.
	As $\lim_{n\to\infty}  \Vert B_w^{p_n+s} x_i\Vert =0$ for $i=1,\dotsc,k$,
	$(x_1,\dotsc,x_k)\in \mathrm{Prox}_k(B_w)\cap U_1\times\dotsb\times U_k$.
	This shows that  $\mathrm{Prox}_k(B_w)$ is dense in $X^k$.

	By $\sup_{m\in\bbz,n\in\bbn}\prod_{j=m}^{m+n} |w_j|=\infty$ and the weight sequence $\{w_j\}$ is bounded, there exist two sequences $\{q_n\}$ in $\bbn$ and $\{m_n\}$ in $\bbz$
	such that $q_n\to\infty$, $|m_n|\to\infty$ as $n\to\infty$ and
	\[
		\lim_{n\to\infty} \prod_{j=m_n}^{m_n+q_n} |w_j|=\infty.
	\]
	We have two cases (1) $m_n>0$ for all $n$ and (2) $m_n<0$ for all $n$.
	We only consider the case (1), as the case (2) is similar.
	Now assume that $m_n>0$ for all $n$. Then
	\[
		\lim_{n\to\infty} \Vert B_w^{q_n} e_{m_n+q_n}\Vert    =\lim_{n\to\infty} \prod_{j=m_n}^{m_n+q_n}|w_j|=\infty.
	\]
	Fix $k\geq 2$ and nonempty open subsets $V_1,\dotsc,V_k$ of $X$.
	Then there exist $s\in \bbn$ and $y_i\in U_i$ such that $y_i$ can be expressed as
	$\sum_{j=-s}^s\beta_j^{(i)} e_j$ for $i=1,\dotsc,k$.
	There exists $\varepsilon>0$ such that $y_i+\alpha e_j\in U_i$ for $i=1,\dotsc,k$, $j>s$ and $|\alpha|<\varepsilon$.
	For each $n\in\bbn$ and $i=1,\dotsc,k$,
	$z_i  = y_i+\frac{i}{k+1} e_{m_n+q_n}$.
	Then
	\[
		\lim_{n\to\infty}\Vert B_w^{q_n} (z_i-z_j) \Vert \geq  \lim_{n\to\infty}\frac{1}{k+1}\Vert B_w^{q_n}e_{m_n+q_n} \Vert =\infty.
	\]
	This show that $(z_1,\dotsc,z_k)\in D_k(B_w)\cap V_1\times \dotsb\times V_k$.
	Hence $D_k(B_w)$ is dense in $X^k$.
	By Theorem~\ref{thm:chaos-equi} $B_w$ is densely uniformly Li-Yorke chaotic.
\end{proof}

According  to Theorem~\ref{thm:bi-backward-LY-chaos} and Lemma~\ref{lem:T-prox-dense-equi}, we have the following consequence about the property of sequences in $\mathbb{K}$, which should be of independent interest.

\begin{coro}
	Let $w=\{w_j\}_{j\in\bbz}$ be a bounded non-zeros sequence in $\mathbb{K}$.
	If
	\[\liminf\limits_{n\to\infty} \prod_{j=-n+1}^0 |w_j|=0, \]
	then there exists a sequences $\{p_n\}$ in  $\bbn$ such that
	\[\lim\limits_{n\to\infty} \prod_{j=-p_n+k+1}^k|w_j|=0\]
	for all $k\in\bbz$.
\end{coro}
\begin{proof}
	We consider the bilateral backward weighted shift $B_w\colon \ell^1(\bbz)\to \ell^1(\bbz)$, $B_w e_i =w_ie_{i-1}$ for all $i\in\bbz$.
	By the proof of (1)$\Rightarrow$(5) of Theorem \ref{thm:bi-backward-LY-chaos},
	For every $k\geq 2$, $\mathrm{Prox}_k(B_w)$ is dense in $\ell^1(\bbz)^k$.
	By Lemma~\ref{lem:T-prox-dense-equi}, there exist a sequence $\{p_n\}$ in  $\bbn$
	and a dense subset $K$ of $\ell^1(\bbz)$ such that for every $x\in K$,
	$\lim_{n\to\infty} B_w^{p_n}x =0$.
	Fix $k\in \bbz$.
	As $K$ is dense in $\ell^1(\bbz)$, there exists $x\in K$ such that $\Vert x-e_k\Vert<1$.
	In particular, $x_k\neq 0$.
	Then
	\[
		\lim_{n\to\infty} \prod_{j=-p_n+k+1}^k |w_j|=\lim_{n\to\infty}   \Vert B_w^{p_n} e_k\Vert \leq \frac{1}{|x_k|} \lim_{n\to\infty} \Vert B_w^{p_n}x\Vert =0.
	\]
This ends the proof.
\end{proof}

In \cite{S2008}*{Theorem 1.1}, Smith constructed an operator $T$ on $\ell^2(\bbn)$ such that the whole space $\ell^2(\bbn)$ is a uniformly Li-Yorke scrambled set.
By linearity, we have the following result: the existence of a ``large'' uniformly Li-Yorke scrambled set implies the whole space is uniformly Li-Yorke scrambled.
Recall that a subset $A$ of $X$ is said to be \emph{locally residual} if there exists a nonempty open subset $U$ of $X$ such that $A\cap U$ is residual in $U$.

\begin{prop}
	Let $T\in\mathcal{B}(X)$.
	If there exists some uniformly Li-Yorke scrambled subset of $X$ which is locally residual, then the whole space $X$ is uniformly Li-Yorke scrambled.
\end{prop}
\begin{proof}
	Let $S$ be a uniformly Li-Yorke scrambled subset of $X$ which is locally residual.
	There exist two sequences $\{p_n\}$ and $\{q_n\}$ in $\bbn$ such that
	for any pair $x,y\in S$ with $x\neq y$, we have
	\[
		\lim_{n\to\infty} \Vert T^{p_n}x -T^{p_n}y  \Vert =0\text{ and }
		\lim_{n\to\infty} \Vert T^{q_n}x -T^{q_n}y\Vert =\infty.
	\]
	As $S$ is locally residual, there exists a vector $z\in X$ and $\delta>0$ such that $S\cap B(z,2\delta)$ is residual in $B(z,2\delta)$.
	Fix a non-zero vector $x\in X$ and pick $\varepsilon>0$ such that $\Vert \varepsilon x\Vert <\delta$.
	By linearity,  $(S+\varepsilon x)\cap B(z+\varepsilon x,2\delta)$ is residual in $B(z+\varepsilon x,2\delta)$.
	Note that $B(z,\delta)\subset B(z+\varepsilon x,2\delta)$,
	then $(S+\varepsilon x)\cap S\cap  B(z,\delta)$ is residual in $B(z,\delta)$.
	Pick a point $y\in (S+\varepsilon x)\cap S$. Then $y,y-\varepsilon x\in S$ and therefore
	\[
		\lim_{n\to\infty} \Vert T^{p_n}x\Vert =0\text{ and }
		\lim_{n\to\infty} \Vert T^{q_n}x\Vert =\infty.
	\]
	This implies that the whole space $X$ is uniformly Li-Yorke scrambled.
\end{proof}

Note that if $A$ is a subset of $X$ of second category with the Baire property then $A$ is locally residual (see e.g.\@
\cite{K1995}*{Proposition I.8.26}), and the collection of sets with the Baire property, the smallest $\sigma$-algebra containing all open sets and
all sets of first category (see e.g.\@ \cite{K1995}*{Proposition I.8.22}). So we have the following consequence.

\begin{coro}
	Let $T\in\mathcal{B}(X)$.
	If there exists some uniformly Li-Yorke scrambled subset of $X$ which is Baire and of second category, then the whole space $X$ is uniformly Li-Yorke scrambled.
\end{coro}

By Corollary~\ref{cor:DULYC-dense-irr-manifold}
dense uniform Li-Yorke chaos implies the existence of a dense irregular manifold.
By Corollary~\ref{coro:unilateral-bws} and Theorem~\ref{thm:bi-backward-LY-chaos}, we know that for unilateral backward weighted shifts and  bilateral backward weighted shifts, Li-Yorke chaos is equivalent to dense uniform Li-Yorke chaos.
But we have the following question.

\begin{question}\label{problem1}
Does there exists an operator $T\in\mathcal{B}(X)$ such that $T$ admits a dense irregular manifold but not densely uniformly Li-Yorke chaotic?
\end{question}

If an operator $T\in\mathcal{B}(X)$ is hypercyclic then $T$ admits a dense irregular manifold.
So a special case of Question~\ref{problem1} is as follows.

\begin{question}
Does the hypercyclicity imply dense uniform Li-Yorke chaos?
\end{question}

\noindent\textbf{Acknowledgments}:
J. Li was partially supported by NSF of China (12222110, 12171298). X. Wang was partially supported by STU Scientific Research Initiation Grant (SRIG, No. NTF22020) and NSF of China (12301230).

\end{document}